\documentclass[12pt]{amsart}
\usepackage{amsfonts}
\usepackage{amsmath,amssymb,amsthm}

\newcommand{\R}{\mathbb R}

\renewcommand{\span}{\mathrm{span}}
\newcommand{\tr}{\mathrm{tr}}

\newtheorem{thm}{Theorem}[section]

\newtheorem{prop}[thm]{Proposition}
\theoremstyle{definition}

\theoremstyle{remark}

\newcommand{\ds}{\displaystyle}

\begin{document}

\title[MARGINALLY TRAPPED SURFACES IN MINKOWSKI SPACE]
{AN INVARIANT THEORY OF MARGINALLY TRAPPED SURFACES IN THE
FOUR-DIMENSIONAL MINKOWSKI SPACE}

\author{Georgi Ganchev and Velichka Milousheva}
\address{Bulgarian Academy of Sciences, Institute of Mathematics and Informatics,
Acad. G. Bonchev Str. bl. 8, 1113 Sofia, Bulgaria}
\email{ganchev@math.bas.bg}
\address{Bulgarian Academy of Sciences, Institute of Mathematics and Informatics,
Acad. G. Bonchev Str. bl. 8, 1113, Sofia, Bulgaria; "L. Karavelov"
Civil Engineering Higher School, 175 Suhodolska Str., 1373 Sofia,
Bulgaria} \email{vmil@math.bas.bg}

\subjclass[2000]{Primary 53A35, Secondary 53B25}
\keywords{Marginally trapped surfaces in the four-dimensional
Minkowski space, lightlike mean curvature vector, Bonnet-type
fundamental theorem, meridian surfaces in Minkowski space}

\begin{abstract}
A marginally trapped surface in the four-dimensional Minkowski space is a spacelike surface whose
mean curvature vector is lightlike at each point.
We associate  a geometrically determined moving frame field to such a surface and using the derivative formulas
for this frame field we obtain seven invariant functions. Our main theorem states that these seven invariants determine
the surface up to a motion in Minkowski space.

We introduce meridian surfaces as  one-parameter systems of meridians of  a rotational hypersurface
in the four-dimensional Minkowski space.
We find all marginally trapped meridian surfaces.
\end{abstract}

\maketitle

\section{Introduction}

The concept of trapped surfaces was introduced by Roger Penrose in \cite{Pen} and it plays an important role in general relativity.
A surface in a 4-dimensional spacetime is called marginally trapped if it is closed, embedded, spacelike and
its mean curvature vector is lightlike at each point of the surface.
These surfaces were defined by Penrose in order to study global properties of spacetime.
In Physics similar or weaker definitions attract attention.
Recently, marginally trapped surfaces have been studied from a mathematical viewpoint. In the mathematical literature,
it is customary to call a codimension-two surface in a 4-dimensional semi-Riemannian manifold  \emph{marginally trapped}
it its mean curvature vector $H$ is lightlike at each point,
and removing the other hypotheses, i.e. the surface does not need to be closed or embedded.
Classification results in 4-dimensional Lorentz manifolds were obtained imposing some extra conditions on the mean curvature vector,
the Gauss curvature or the second fundamental form.
Marginally trapped surfaces with positive relative nullity in Lorenz space forms were classified in \cite{Chen-Veken-1}.
The non-existence of marginally trapped surfaces in Robertson-Walker spaces with positive relative nullity was shown in
\cite{Chen-Veken-2}.
Marginally trapped surfaces with parallel mean curvature vector in Lorenz space forms were classified in \cite{Chen-Veken-3}.
In \cite{Haesen-Ort-1} marginally trapped surfaces which are invariant under a boost transformation
in 4-dimensional Minkowski space were studied,
and marginally trapped surfaces in Minkowski 4-space which are invariant under spacelike rotations were classified in \cite{Haesen-Ort-2}.
The classification of marginally trapped surfaces in Minkowski 4-space which are invariant under a group of screw rotations
(a group of Lorenz rotations with an invariant lightlike direction) is obtained in \cite{Haesen-Ort-3}.

In this paper, we consider marginally trapped surfaces in the four-dimensional
Minkowski space $\R^4_1$. Our study is based on the geometrically introduced  invariant
linear map of Weingarten-type in the tangent plane at any point of
the surface under consideration. This allows us to introduce
principal lines and geometrically
determined invariant moving frame field. Writing
derivative formulas of Frenet-type for this frame field, we obtain
seven invariant functions and prove a fundamental theorem of
Bonnet-type, stating that these seven invariants under some
natural conditions determine the surface up to a motion in $\R^4_1$.

We apply our theory to spacelike surfaces lying on rotational hypersurfaces in $\R^4_1$.
Considering rotational hypersurfaces  with timelike or  spacelike axis
we construct special two-dimensional  surfaces which are one-parameter systems of meridians
of the rotational hypersurface (meridian surfaces).
We find all meridian surfaces which are marginally trapped.

\section{Preliminaries}

In \cite{GM5} we considered the local theory of spacelike surfaces
in the four-dimensional Minkowski space $\R^4_1$. The basic
feature of our treatment of these surfaces was the introduction of an
invariant linear map of Weingarten-type in the tangent plane at
any point of the surface, following the  approach to the theory of surfaces in $\R^4$ \cite{GM1,GM4}.
 Studying surfaces in the Euclidean space
$\R^4$, in \cite{GM1} we introduced a linear map $\gamma$ of
Weingarten-type, which plays a similar role in the theory of
surfaces in $\R^4$ as the Weingarten map in the theory of surfaces
in $\R^3$. The map $\gamma$ generates the corresponding second
fundamental form $II$ at any point of the surface in the standard
way. We gave a geometric interpretation of the second
fundamental form and the Weingarten map of the surface in
\cite{GM3}.

\vskip 2mm
Let  $\R^4_1$ be the Minkowski space endowed with the metric
$\langle , \rangle$ of signature $(3,1)$ and $Oe_1e_2e_3e_4$ be a
fixed orthonormal coordinate system in $\R^4_1$, i.e. $e_1^2 =
e_2^2 = e_3^2 = 1, \, e_4^2 = -1$, giving the orientation of
$\R^4_1$. The standard flat metric is given in local coordinates by
$dx_1^2 + dx_2^2 + dx_3^2 -dx_4^2.$

A surface $M^2$ in $\R^4_1$ is said to be
\emph{spacelike} if $\langle , \rangle$ induces  a Riemannian
metric $g$ on $M^2$. Thus at each point $p$ of a spacelike surface
$M^2$ we have the following decomposition
$$\R^4_1 = T_pM^2 \oplus N_pM^2$$
with the property that the restriction of the metric
$\langle , \rangle$ onto the tangent space $T_pM^2$ is of
signature $(2,0)$, and the restriction of the metric $\langle ,
\rangle$ onto the normal space $N_pM^2$ is of signature $(1,1)$.

Denote by $\nabla'$ and $\nabla$ the Levi Civita connections on $\R^4_1$ and $M^2$, respectively.
Let $x$ and $y$ denote vector fields tangent to $M$ and let $\xi$ be a normal vector field.
Then the formulas of Gauss and Weingarten give a decomposition of the vector fields $\nabla'_xy$ and
$\nabla'_x \xi$ into a tangent and a normal component:
$$\begin{array}{l}
\vspace{2mm}
\nabla'_xy = \nabla_xy + \sigma(x,y);\\
\vspace{2mm}
\nabla'_x \xi = - A_{\xi} x + D_x \xi,
\end{array}$$
which define the second fundamental tensor $\sigma$, the normal connection $D$ and the shape operator $A_{\xi}$ with respect to $\xi$.
The mean curvature vector  field $H$ of the surface $M^2$ is defined as $H = \ds{\frac{1}{2}\,  \tr\, \sigma}$,
i.e. given a local orthonormal frame $\{x,y\}$ of the tangent bundle, $H = \ds{\frac{1}{2} \left(\sigma(x,x) +  \sigma(y,y)\right)}$.

Let
$M^2: z=z(u,v), \,\, (u,v) \in \mathcal{D}$ $(\mathcal{D} \subset \R^2)$  be a local parametrization on a
spacelike surface in $\R^4_1$.
The tangent space at an arbitrary point $p=z(u,v)$ of $M^2$ is $T_pM^2 = \span \{z_u,z_v\}$. Since $M^2$ is spacelike,
$\langle z_u,z_u \rangle > 0$, $\langle z_v,z_v \rangle > 0$.
We use the standard denotations
$E(u,v)=\langle z_u,z_u \rangle, \; F(u,v)=\langle z_u,z_v
\rangle, \; G(u,v)=\langle z_v,z_v \rangle$ for the coefficients
of the first fundamental form
$$I(\lambda,\mu)= E \lambda^2 + 2F \lambda \mu + G \mu^2,\,\,
\lambda, \mu \in \R.$$
 Since $I(\lambda, \mu)$ is positive
definite we set $W=\sqrt{EG-F^2}$.
We choose a normal frame field $\{n_1, n_2\}$ such that $\langle
n_1, n_1 \rangle =1$, $\langle n_2, n_2 \rangle = -1$, and the
quadruple $\{z_u,z_v, n_1, n_2\}$ is positively oriented in
$\R^4_1$.
Then we have the following derivative formulas:
$$\begin{array}{l}
\vspace{2mm} \nabla'_{z_u}z_u=z_{uu} = \Gamma_{11}^1 \, z_u +
\Gamma_{11}^2 \, z_v + c_{11}^1\, n_1 - c_{11}^2\, n_2;\\
\vspace{2mm} \nabla'_{z_u}z_v=z_{uv} = \Gamma_{12}^1 \, z_u +
\Gamma_{12}^2 \, z_v + c_{12}^1\, n_1 - c_{12}^2\, n_2;\\
\vspace{2mm} \nabla'_{z_v}z_v=z_{vv} = \Gamma_{22}^1 \, z_u +
\Gamma_{22}^2 \, z_v + c_{22}^1\, n_1 - c_{22}^2\, n_2,\\
\end{array}$$
where $\Gamma_{ij}^k$ are the Christoffel's symbols and the functions $c_{ij}^k, \,\, i,j,k = 1,2$ are given by
$$\begin{array}{lll}
\vspace{2mm}
c_{11}^1 = \langle z_{uu}, n_1 \rangle; & \qquad   c_{12}^1 = \langle z_{uv}, n_1 \rangle; & \qquad  c_{22}^1 = \langle z_{vv}, n_1 \rangle;\\
\vspace{2mm}
c_{11}^2 = \langle z_{uu}, n_2 \rangle; & \qquad  c_{12}^2 = \langle z_{uv}, n_2 \rangle; & \qquad
c_{22}^2 = \langle z_{vv}, n_2 \rangle.
\end{array} $$

Obviously, the surface $M^2$ lies in a 2-plane if and only if
$M^2$ is totally geodesic, i.e. $c_{ij}^k=0, \; i,j,k = 1, 2.$ So,
we assume that at least one of the coefficients $c_{ij}^k$ is not
zero.

\vskip 2mm
Considering the tangent
space $T_pM^2$ at a point $p \in M^2$, in \cite{GM5} we introduced an invariant
$\zeta_{\,g_1,g_2}$ of a pair of two tangents $g_1$, $g_2$ using
the second fundamental tensor $\sigma$ of $M^2$. By means of this
invariant we defined conjugate, asymptotic, and principal
tangents.
The second fundamental form $II$ of the surface $M^2$ at a point
$p \in M^2$ is introduced on the base of conjugacy of two tangents
at the point.
The coefficients $L, M, N$ of the second fundamental form $II$ are determined as follows:
\begin{equation} \label{Eq-1}
L = \ds{\frac{2}{W}} \left|%
\begin{array}{cc}
\vspace{2mm}
  c_{11}^1 & c_{12}^1 \\
  c_{11}^2 & c_{12}^2 \\
\end{array}%
\right|; \quad
M = \ds{\frac{1}{W}} \left|%
\begin{array}{cc}
\vspace{2mm}
  c_{11}^1 & c_{22}^1 \\
  c_{11}^2 & c_{22}^2 \\
\end{array}%
\right|; \quad
N = \ds{\frac{2}{W}} \left|%
\begin{array}{cc}
\vspace{2mm}
  c_{12}^1 & c_{22}^1 \\
  c_{12}^2 & c_{22}^2 \\
\end{array}%
\right|.
\end{equation}

The second fundamental form $II$ determines an invariant linear
map $\gamma$  of Weingarten-type at any point of the surface,
which generates  two invariant functions:
$$k := \det \gamma = \frac{LN - M^2}{EG - F^2}, \qquad
\varkappa := -\frac{1}{2}\,{\rm tr}\, \gamma =
\frac{EN+GL-2FM}{2(EG-F^2)}.$$
The functions $k$ and $\varkappa$
are invariant under changes of the parameters of the surface and
changes of the normal frame field.
 The sign of $k$ is invariant under congruences
 and the sign of $\varkappa$ is invariant under motions
in $\R^4_1$. However, the sign of $\varkappa$ changes under
symmetries with respect to a hyperplane in $\R^4_1$. We proved that
the invariant $\varkappa$ is the curvature of the normal
connection of the surface. The number of asymptotic tangents at a
point of $M^2$ is determined by the sign of the invariant $k$. In
the case $k=0$ there exists a one-parameter family of asymptotic
lines, which are principal.

It is interesting to note that the
''umbilical'' points, i.e. points at which the coefficients of the
first and the second fundamental forms are proportional, are
exactly the points at which the mean curvature vector $H$ is zero.
So, the spacelike surfaces consisting of ''umbilical'' points in
$\R^4_1$ are exactly the minimal surfaces. Minimal spacelike
surfaces are characterized in terms of the invariants $k$ and
$\varkappa$ by the equality $\varkappa^2 - k =0$.

Analogously to $\R^3$ and $\R^4$, the invariants $k$ and
$\varkappa$ divide the points of $M^2$ into four types: flat,
elliptic, hyperbolic and parabolic points. The surfaces consisting
of flat points  are characterized by the conditions $k= \varkappa
= 0$, or equivalently $L=M=N=0$. We gave a local geometric description of spacelike surfaces
 consisting of flat points whose mean curvature vector
at any point is a non-zero spacelike vector or  timelike vector,
proving that any such a surface either lies in a hyperplane of
$\R^4_1$ or is part of a developable ruled surface in $\R^4_1$ \cite{GM5}.

Using the introduced principal lines on a spacelike surface in
$\R^4_1$ whose mean curvature vector at any point is a non-zero
spacelike vector or  timelike vector, we found a geometrically
determined moving frame field  on such a surface. Writing the
derivative formulas of Frenet-type for this frame field, we
obtained eight invariant functions and proved a fundamental
theorem of Bonnet-type, stating that these eight invariants under
some natural conditions determine the surface up to a motion in
$\R^4_1$.

In the present paper we shall apply the same idea for
developing the invariant theory of
spacelike surfaces   in
$\R^4_1$ whose mean curvature vector at any point is a lightlike vector, i.e.  marginally trapped surfaces.

\section{Invariants of a marginally trapped surface} \label{S:Invariants}

Let $M^2: z=z(u,v), \,\, (u,v) \in \mathcal{D}$  be a marginally trapped surface.
Then the mean curvature vector is lightlike at each point of the surface, i.e. $\langle H, H \rangle = 0$.
Thus there exists a pseudo-orthonormal normal frame field $\{n_1, n_2\}$, such that $n_1 = H$ and
$$\langle n_1, n_1 \rangle = 0; \quad \langle n_2, n_2 \rangle = 0; \quad \langle n_1, n_2 \rangle = -1.$$
We assume that $M^2$  is  free of flat points, i.e. $(L,M,N) \neq (0,0,0)$.
Then at each point of the surface there exist principal lines and
without loss of generality we assume that $M^2$ is parameterized by principal lines. Let us denote $x = \ds{\frac{z_u}{\sqrt{E}}}$,
$y = \ds{\frac{z_v}{\sqrt{G}}}$.
Then $\sigma(x,x)$ and $\sigma(y,y)$ are collinear with the mean curvature vector field.
So, we have the following decompositions
\begin{equation}\label{Eq-2}
\begin{array}{l}
\vspace{2mm}
\sigma(x,x) = (1+ \nu) \,n_1; \\
\vspace{2mm}
\sigma(x,y) = \lambda\,n_1 + \mu\,n_2;  \\
\vspace{2mm}
 \sigma(y,y) = (1- \nu) \,n_1,
\end{array}
\end{equation}
where $\nu, \lambda, \mu$ are invariant functions, determined by the principal directions $x, y$, and  the mean curvature vector field
$n_1 = H$ as follows:
$$\nu = - \langle \frac{\sigma(x,x) - \sigma(y,y)}{2}, n_2 \rangle, \quad
\lambda =-  \langle \sigma(x,y), n_2\rangle, \quad \mu = -\langle \sigma(x,y), n_1\rangle.$$

Using \eqref{Eq-1} and \eqref{Eq-2} we calculate the coefficients $L, M, N$ of the second fundamental form and find the invariants
$k$, $\varkappa$ and the Gauss curvature $K$  of $M^2$ expressed by the functions $\nu, \lambda$, and $\mu$:

\begin{equation} \label{Eq-3}
k = 4 \mu^2 (\nu^2 -1); \qquad \varkappa = - 2 \mu \nu; \qquad K = 2 \lambda \mu.
\end{equation}
Since $H \neq 0$, we have $\varkappa^2 - k > 0$.
From \eqref{Eq-3} we get the relations:
$$
4 \mu^2 = \varkappa^2 - k; \qquad \lambda = \ds{\frac{K}{\sqrt{\varkappa^2 - k}}},
\qquad \nu = \ds{- \frac{\varkappa}{\sqrt{\varkappa^2 - k}}}.$$
The last equalities show the geometric meaning of the invariants $\nu, \lambda$, and $\mu$.
Note that $\mu \neq 0$, since $\varkappa^2 - k > 0$.

\vskip 2mm
Recall that a surface is called flat if the Gauss curvature $K$ vanishes;
$M^2$ is with flat normal connection if and only if $\varkappa = 0$;
 $M^2$ consists of parabolic points if and only if $k = 0$. So, equalities \eqref{Eq-3}
imply the following statements.

\begin{prop}\label{P:Flat surface}
Let $M^2$ be a marginally trapped surface free of flat points.
Then $M^2$ is a flat surface if and only if  $\lambda = 0$.
\end{prop}

\begin{prop}\label{P:Flat normal connection}
Let $M^2$ be a marginally trapped surface free of flat points.
Then $M^2$ is a surface with flat normal connection if and only if  $\nu = 0$.
\end{prop}

\begin{prop}\label{P:Parabolic points}
Let $M^2$ be a marginally trapped surface free of flat points.
Then $M^2$ is a surface consisting of parabolic points if and only if $\nu = \pm 1$.
\end{prop}

The flat marginally trapped surfaces can also be characterized
in terms of the so called null allied mean curvature vector field.
The \emph{allied vector field} of a  normal vector field $\xi$ of an $n$-dimensional submanifold
 $M^n$ of $(n+m)$-dimensional Riemannian manifold $\widetilde{M}^{n+m}$ is defined by B.-Y. Chen  \cite{Chen1} by
the formula
$$a(\xi) = \ds{\frac{\|\xi\|}{n} \sum_{k=2}^m \{\tr(A_1 \circ A_k)\}\xi_k},$$
where $\{\xi_1 = \ds{\frac{\xi}{\|\xi\|}},\xi_2,\dots, \xi_m \}$
is a local orthonormal frame of the normal bundle of $M^n$, and $A_i =
A_{\xi_i}, \,\, i = 1,\dots, m$ is the shape operator with respect
to $\xi_i$. In particular, the allied vector field $a(H)$ of the
mean curvature vector field $H$ is a well-defined normal vector
field which is  called the \emph{allied
mean curvature vector field} of $M^n$ in $\widetilde{M}^{n+m}$.
B.-Y. Chen defined  the $\mathcal{A}$-submanifolds to be those
submanifolds of $\widetilde{M}^{n+m}$ for which
 $a(H)$ vanishes identically \cite{Chen1}.
In \cite{GVV1,GVV2} the $\mathcal{A}$-submanifolds are called
\emph{Chen submanifolds}. It is easy to see that minimal
submanifolds, pseudo-umbilical submanifolds and hypersurfaces are
Chen submanifolds. These Chen submanifolds are said to be \emph{trivial
$\mathcal{A}$-submanifolds}.

In \cite{Haesen-Ort-3} S. Haesen and M. Ortega extended the notion of allied mean curvature vector field to the case
when the normal space is a two-dimensional Lorenz space and the mean curvature vector field is lightlike as follows.
Denote by $\{H, H^{\bot}\}$ a pseudo-orthonormal basis of the normal space such that
$\langle H, H \rangle = 0; \,\, \langle H^{\bot}, H^{\bot} \rangle = 0; \,\, \langle H, H^{\bot} \rangle = -1.$
The \emph{null allied mean curvature vector field} is defined as
\begin{equation} \label{Eq-4}
a(H) = \ds{\frac{1}{2}\, \tr(A_H \circ A_{H^{\bot}}) \, H^{\bot}}.
\end{equation}

Now, if $M^2$ is a marginally trapped surface, then using equalities \eqref{Eq-2} we get
$$A_H = A_{n_1} = \left(%
\begin{array}{cc}
  0 & -\mu \\
  -\mu & 0 \\
\end{array}%
\right); \qquad
A_{H^{\bot}} = A_{n_2} = \left(%
\begin{array}{cc}
  -(1+\nu) & -\lambda \\
  -\lambda & -(1-\nu) \\
\end{array}%
\right).
$$
Applying formula \eqref{Eq-4} and using \eqref{Eq-3},
we obtain that the null allied mean curvature vector field of $M^2$ is expressed as follows:
$$a(H) = \lambda \mu \, n_2 = \ds{\frac{K}{2} \, n_2}.$$
Thus we obtain the following result.

\begin{prop}\label{P:Null allied mean CVF}
Let $M^2$ be a marginally trapped surface free of flat points.
Then $M^2$ is a flat surface if and only if  $M^2$ has  vanishing null allied mean curvature vector field.
\end{prop}

In \cite{Haesen-Ort-3} the same result is proved for screw invariant marginally trapped surfaces in $\R^4_1$.

\section{Fundamental theorem}\label{S:Fundamental Thm}

In the local theory of surfaces in Euclidean space  a statement of
significant importance is a theorem of Bonnet-type giving the
natural conditions under which the surface is determined up to a
motion. A theorem of this type was proved for surfaces with flat
normal connection by B.-Y. Chen in \cite{Chen1}. In \cite{GM4} we
proved a fundamental theorem of Bonnet-type for surfaces
in $\R^4$ free of minimal points. In \cite{GM5} we considered spacelike
surfaces in $\R^4_1$ whose mean curvature vector at any point is a
non-zero spacelike vector or  timelike vector. Using the
geometric moving frame field of Frenet-type on such a
surface  and the corresponding derivative formulas, we proved a fundamental
theorem of Bonnet-type for this class of surfaces, stating that any such a surface is determined up to a motion in
$\R^4_1$ by eight invariant functions satisfying some natural conditions.

In this section we shall consider marginally trapped surfaces, i.e. spacelike surfaces
whose mean curvature vector at any point is a lightlike
vector. Let $M^2$  be such a surface. We assume that $M^2$  is  free of flat points, i.e. $(L,M,N) \neq (0,0,0)$,
and $M^2$ is parameterized by principal lines. Let $\{x,y\}$ be the principal tangent vector fields and
$\{n_1, n_2\}$ be the pseudo-orthonormal normal frame field, defined in Section \ref{S:Invariants}.
Thus we obtain
a special frame field $\{x,y,n_1,n_2\}$
at each point $p \in M^2$, such that $x, y$ are unit spacelike vector fields collinear with the principal directions;
$n_1, n_2$ are lightlike vectors, $\langle n_1,\, n_2 \rangle = -1$, and $n_1$ is the mean curvature vector field.
We call such a frame field a \emph{geometric frame field} of $M^2$.

With respect to this frame
field we have the following Frenet-type derivative
formulas of $M^2$:
\begin{equation} \label{Eq-5}
\begin{array}{ll}
\vspace{2mm} \nabla'_xx=\quad \quad \quad \gamma_1\,y+\,(1+\nu)\,n_1;
& \qquad
\nabla'_x n_1= \quad\quad \quad \quad  \mu\,y+\beta_1\,n_1;\\
\vspace{2mm} \nabla'_xy=-\gamma_1\,x\quad \quad \; + \,\quad  \lambda\,n_1
\; + \mu\,n_2;  & \qquad
\nabla'_y n_1=\mu \,x \quad\quad \quad \quad +\beta_2\,n_1;\\
\vspace{2mm} \nabla'_yx=\quad\quad \;\, -\gamma_2\,y  + \quad \lambda\,n_1
\; +\mu\,n_2;  & \qquad
\nabla'_xn_2= (1+\nu) \, x + \lambda \,y \quad \quad  -\beta_1\,n_2;\\
\vspace{2mm} \nabla'_yy=\;\;\gamma_2\,x \quad\quad\;\;\,
+(1-\nu)\,n_1; & \qquad \nabla'_y n_2= \lambda \, x  + (1-\nu) \,y
 \quad \quad -\beta_2\,n_2,
\end{array}
\end{equation}
where  $\gamma_1 = - y(\ln \sqrt{E}), \,\, \gamma_2 = - x(\ln
\sqrt{G})$, $\beta_1 = - \langle \nabla'_x n_1, n_2\rangle$, $\beta_2 = -
\langle \nabla'_y n_1, n_2\rangle$, and $\mu \neq 0$.

Using that $R'(x,y,n_1) = 0$, $R'(x,y,n_2) = 0$, and $R'(x,y,x) = 0$,
from \eqref{Eq-5} we get the following integrability conditions:
$$\begin{array}{l}
\vspace{2mm}
 x(\mu) - 2\mu\, \gamma_2  - \mu\,\beta_1 =0;\\
\vspace{2mm}
y(\mu) - 2\mu\, \gamma_1 - \mu\,\beta_2 = 0;\\
\vspace{2mm}
 x(\gamma_2) +
y(\gamma_1) - \left((\gamma_1)^2
+ (\gamma_2)^2\right) - 2\lambda \, \mu = 0;\\
\vspace{2mm}
x(\lambda) - y(\nu) - 2\lambda\, \gamma_2 + 2\nu\, \gamma_1 + \lambda\,\beta_1 - (1+ \nu)\,\beta_2 = 0;\\
\vspace{2mm}
x(\nu) + y(\lambda) - 2\lambda\, \gamma_1 - 2\nu\, \gamma_2 - (1- \nu)\,\beta_1 + \lambda\,\beta_2 = 0;\\
\vspace{2mm}
x(\beta_2) - y(\beta_1) + 2 \nu\,\mu + \gamma_1\,\beta_1 - \gamma_2\,\beta_2 = 0.
\end{array}$$

Having in mind that $x = \displaystyle{\frac{z_u}{\sqrt{E}}, \, y
= \frac{z_v}{\sqrt{G}}}$, we can rewrite the above equalities as follows:
$$\begin{array}{l}
\vspace{2mm} 2\mu\, \gamma_2 + \mu\,\beta_1 =
\displaystyle{\frac{1}{\sqrt{E}} \, \mu_u};\\
\vspace{2mm} 2\mu\, \gamma_1+ \mu\,\beta_2 =
\displaystyle{\frac{1}{\sqrt{G}}\, \mu_v};\\
\vspace{2mm}
2\lambda\, \mu =
\displaystyle{\frac{1}{\sqrt{E}}\,(\gamma_2)_u
+ \frac{1}{\sqrt{G}}\,(\gamma_1)_v - \left((\gamma_1)^2 + (\gamma_2)^2\right)};\\
\vspace{2mm} 2\lambda\, \gamma_2 - 2\nu\, \gamma_1 - \lambda\,\beta_1 + (1+ \nu)\,\beta_2 =
\displaystyle{\frac{1}{\sqrt{E}}\, \lambda_u - \frac{1}{\sqrt{G}}\,\nu_v};\\
\vspace{2mm} 2\lambda\, \gamma_1 + 2\nu\, \gamma_2 + (1 - \nu)\,\beta_1 - \lambda\,\beta_2 =
\displaystyle{\frac{1}{\sqrt{E}}\,\nu_u + \frac{1}{\sqrt{G}}\,\lambda_v};\\
\gamma_1\,\beta_1 - \gamma_2\,\beta_2 + 2 \nu\,\mu  =
\displaystyle{ - \frac{1}{\sqrt{E}}\,(\beta_2)_u +
\frac{1}{\sqrt{G}}\,(\beta_1)_v}.
\end{array}$$

The condition $\mu_u \,\mu_v \neq 0$ is equivalent to $(2\, \gamma_2 + \beta_1) (2\, \gamma_1 + \beta_2) \neq 0$.
So, if $\mu_u \,\mu_v \neq 0$, then
$$\sqrt{E} = \displaystyle{\frac{\mu_u}
{\mu (2\, \gamma_2 + \beta_1)}}; \quad
\sqrt{G} = \displaystyle{\frac{\mu_v}{\mu (2\, \gamma_1 + \beta_2)}}.$$

\vskip 3mm We shall prove the following  Bonnet-type theorem for marginally trapped
 surfaces in $\R^4_1$ free of flat points.

\begin{thm}\label{T:Main Theorem}
Let $\gamma_1, \, \gamma_2, \, \nu, \, \lambda, \, \mu,
\, \beta_1, \beta_2$ be smooth functions, defined in a domain
$\mathcal{D}, \,\, \mathcal{D} \subset {\R}^2$, and satisfying the
conditions
\begin{equation} \label{Eq-6}
\begin{array}{l}
\vspace{2mm}
 \displaystyle{\frac{\mu_u}{\mu (2\, \gamma_2 + \beta_1)}}>0; \qquad \qquad \qquad
 \displaystyle{\frac{\mu_v}{\mu (2\, \gamma_1 + \beta_2)}}>0;\\
\vspace{2mm}
- \gamma_1 \sqrt{E} \sqrt{G} = (\sqrt{E})_v; \qquad \qquad
- \gamma_2 \sqrt{E} \sqrt{G} = (\sqrt{G})_u;\\
\vspace{2mm} 2\lambda\, \mu =
\displaystyle{\frac{1}{\sqrt{E}}\,(\gamma_2)_u
+ \frac{1}{\sqrt{G}}\,(\gamma_1)_v - \left((\gamma_1)^2 + (\gamma_2)^2\right)};\\
\vspace{2mm} 2\lambda\, \gamma_2 - 2\nu\, \gamma_1 - \lambda\,\beta_1 + (1+ \nu)\,\beta_2 =
\displaystyle{\frac{1}{\sqrt{E}}\, \lambda_u - \frac{1}{\sqrt{G}}\,\nu_v};\\
\vspace{2mm} 2\lambda\, \gamma_1 + 2\nu\, \gamma_2 + (1 - \nu)\,\beta_1 - \lambda\,\beta_2 =
\displaystyle{\frac{1}{\sqrt{E}}\,\nu_u + \frac{1}{\sqrt{G}}\,\lambda_v};\\
\gamma_1\,\beta_1 - \gamma_2\,\beta_2 + 2 \nu\,\mu  =
\displaystyle{ - \frac{1}{\sqrt{E}}\,(\beta_2)_u +
\frac{1}{\sqrt{G}}\,(\beta_1)_v},
\end{array}
\end{equation}
where $\sqrt{E} = \ds{\frac{\mu_u} {\mu (2\, \gamma_2 +
\beta_1)}}$, $\sqrt{G} = \ds{\frac{\mu_v}{\mu (2\, \gamma_1 +
\beta_2)}}$. Let $\{x_0, \, y_0, \, (n_1)_0,\, (n_2)_0\}$ be
vectors  at a point $p_0 \in \R^4_1$, such that $x_0$, $y_0$ are
unit spacelike vectors, $\langle x_0, y_0 \rangle = 0$, $(n_1)_0,
(n_2)_0$ are lightlike vectors, and $\langle (n_1)_0, (n_2)_0
\rangle = -1$. Then there exist a subdomain ${\mathcal{D}}_0
\subset \mathcal{D}$ and a unique marginally trapped surface $M^2:
z = z(u,v), \,\, (u,v) \in {\mathcal{D}}_0$ free of flat points,
such that $M^2$ passes through $p_0$, the functions
 $\gamma_1, \, \gamma_2, \, \nu, \,
\lambda, \, \mu, \, \beta_1, \beta_2$ are the geometric functions
of $M^2$ and $\{x_0, \, y_0, \, (n_1)_0,\, (n_2)_0\}$ is the geometric
frame of $M^2$ at the point $p_0$.
\end{thm}

\vskip 2mm \noindent \emph{Proof:} We consider the following
system of partial differential equations for the unknown vector
functions $x = x(u,v), \, y = y(u,v), \,n_1 = n_1(u,v), \,n_2 = n_2(u,v)$
in $\R^4_1$:
\begin{equation} \label{E:Eq-7}
\begin{array}{ll}
\vspace{2mm} x_u = \sqrt{E} \left(\gamma_1\, y + (1+\nu)\,n_1\right)
& \quad x_v = \sqrt{G}\left(- \gamma_2\, y + \lambda\, n_1 + \mu\, n_2\right)\\
\vspace{2mm} y_u = \sqrt{E} \left(- \gamma_1\, x + \lambda\, n_1 + \mu\, n_2\right)  &
\quad y_v = \sqrt{G}\left(\gamma_2\, x +(1- \nu)\, n_1 \right)\\
\vspace{2mm} (n_1)_u =  \sqrt{E}\left(\mu\, y + \beta_1\, n_1\right)  &
\quad (n_1)_v =  \sqrt{G}\left(\mu\, x + \beta_2\, n_1 \right)\\
\vspace{2mm} (n_2)_u = \sqrt{E}\left((1+\nu)\, x +  \lambda\, y - \beta_1\, n_2\right)&
\quad (n_2)_v = \sqrt{G}\left(\lambda\, x +   (1- \nu)\,y - \beta_2\, n_2 \right)
\end{array}
\end{equation}
We denote
$$Z =
\left(%
\begin{array}{c}
  x \\
  y \\
  n_1 \\
  n_2 \\
\end{array}%
\right); \quad
A = \sqrt{E} \left(%
\begin{array}{cccc}
  0 & \gamma_1 & (1+\nu) & 0 \\
  -\gamma_1 & 0 & \lambda &  \mu \\
  0 & \mu & \beta_1 & 0 \\
  (1+\nu) & \lambda & 0 & -\beta_1 \\
\end{array}%
\right);$$
$$B = \sqrt{G}
\left(%
\begin{array}{cccc}
  0 & -\gamma_2 & \lambda &  \mu \\
  \gamma_2 & 0 & (1-\nu) & 0 \\
  \mu & 0 & \beta_2 & 0 \\
  \lambda & (1-\nu) & 0 & -\beta_2 \\
\end{array}%
\right).$$ Then system \eqref{E:Eq-7} can be rewritten in the form:
\begin{equation} \label{E:Eq-8}
\begin{array}{l}
\vspace{2mm}
Z_u = A\,Z,\\
\vspace{2mm} Z_v = B\,Z.
\end{array}
\end{equation}
The integrability conditions of  \eqref{E:Eq-8}  are
$$Z_{uv} = Z_{vu},$$
i.e.
\begin{equation}\label{E:Eq-9}
\displaystyle{\frac{\partial a_i^k}{\partial v} - \frac{\partial b_i^k}{\partial u}
+ \sum_{j=1}^{4}(a_i^j\,b_j^k - b_i^j\,a_j^k) = 0, \quad i,k = 1,
\dots, 4,}
\end{equation}
 where $a_i^j$ and $b_i^j$ are the
elements of the matrices $A$ and $B$. Using \eqref{Eq-6}  we obtain that
 equalities \eqref{E:Eq-9}  are fulfilled. Hence, there exist a subset
$\mathcal{D}_1 \subset \mathcal{D}$ and unique vector functions $x
= x(u,v), \, y = y(u,v), \,n_1 = n_1(u,v), \,n_2 = n_2(u,v), \,\, (u,v)
\in \mathcal{D}_1$, which satisfy system \eqref{E:Eq-7} and the conditions
$$x(u_0,v_0) = x_0, \quad y(u_0,v_0) = y_0, \quad n_1(u_0,v_0) = (n_1)_0, \quad n_2(u_0,v_0) = (n_2)_0.$$

We shall prove that for each $(u,v) \in \mathcal{D}_1$ the vectors
$x(u,v), \, y(u,v), \,n_1(u,v), \,n_2(u,v)$ have the following properties:
$x(u,v)$, $y(u,v)$ are unit spacelike vectors, $\langle x(u,v), y(u,v) \rangle =0$, $n_1(u,v),\, n_2(u,v)$
are lightlike vectors, and
$\langle n_1,\, n_2 \rangle = -1$. Let us consider the following
functions:
$$\begin{array}{lll}
\vspace{2mm}
  \varphi_1 = \langle x,x \rangle - 1; & \qquad \varphi_5 =
  \langle x,y \rangle; & \qquad \varphi_8 = \langle y,n_1 \rangle; \\
\vspace{2mm}
  \varphi_2 = \langle y, y \rangle - 1; & \qquad \varphi_6 =
  \langle x,n_1 \rangle; & \qquad \varphi_9 = \langle y,n_2 \rangle; \\
\vspace{2mm}
  \varphi_3 = \langle n_1, n_1 \rangle; & \qquad \varphi_7 =
  \langle x,n_2 \rangle; & \qquad \varphi_{10} = \langle n_1,n_2 \rangle + 1; \\
\vspace{2mm}
  \varphi_4 = \langle n_2, n_2 \rangle; &   &  \\
\end{array}$$
defined for each $(u,v) \in \mathcal{D}_1$. Using that $x(u,v), \,
y(u,v), \,n_1(u,v), \,n_2(u,v)$ satisfy \eqref{E:Eq-7}, we obtain  the system
\begin{equation} \label{E:Eq-10}
\begin{array}{lll}
\vspace{2mm}
\displaystyle{\frac{\partial \varphi_i}{\partial u} = \alpha_i^j \, \varphi_j},\\
\vspace{2mm} \displaystyle{\frac{\partial \varphi_i}{\partial v} =
\beta_i^j \, \varphi_j};
\end{array} \qquad i = 1, \dots, 10,
\end{equation}
where $\alpha_i^j, \beta_i^j, \,\, i,j = 1, \dots, 10$ are
functions of $(u,v) \in \mathcal{D}_1$. System \eqref{E:Eq-10} is a linear
system of partial differential equations for the functions
$\varphi_i(u,v), \,\,i = 1, \dots, 10, \,\,(u,v) \in
\mathcal{D}_1$, satisfying $\varphi_i(u_0,v_0) = 0, \,\,i = 1,
\dots, 10$. Hence, $\varphi_i(u,v) = 0, \,\,i = 1, \dots, 10$ for
each $(u,v) \in \mathcal{D}_1$. Consequently, the quadruple
 $\{x(u,v), \, y(u,v), \,n_1(u,v), \,n_2(u,v)\}$ has the properties mentioned above.

Now, let us consider the  system
\begin{equation}\label{E:Eq-11}
\begin{array}{lll}
\vspace{2mm}
z_u = \sqrt{E}\, x\\
\vspace{2mm} z_v = \sqrt{G}\, y
\end{array}
\end{equation}
of partial differential equations for the vector function
$z(u,v)$. Using \eqref{Eq-6} and \eqref{E:Eq-7} we get that the integrability
conditions $z_{uv} = z_{vu}$ of system \eqref{E:Eq-11}
 are fulfilled. Hence,  there exist a subset $\mathcal{D}_0 \subset \mathcal{D}_1$ and
a unique vector function $z = z(u,v)$, defined for $(u,v) \in
\mathcal{D}_0$ and satisfying $z(u_0, v_0) = p_0$.

Consequently, the surface $M^2: z = z(u,v), \,\, (u,v) \in
\mathcal{D}_0$ satisfies the assertion of the theorem. \qed

\vskip 3mm
Marginally trapped surfaces for which $\beta_1 = \beta_2 =0$ have parallel mean curvature vector field,
i.e. $DH =0$ holds identically.
The class of marginally trapped surfaces with parallel mean curvature vector field,
was classified by B.-Y. Chen and J. Van Der Veken in \cite{Chen-Veken-3}.

\vskip 5mm
\section{ Meridian surfaces in $\R^4_1$}\label{S:Examples}

In \cite{GM4} we constructed a family of surfaces lying on a standard rotational hypersurface
in the four-dimensional Euclidean space
$\R^4$. These surfaces are one-parameter systems of meridians
of the rotational hypersurface, that is why we called them \emph{meridian surfaces}.
We described the meridian surfaces with constant Gauss curvature, with
constant mean curvature, and with constant invariant $k$.

In this section we shall use the same idea to construct a special family of two-dimensional  spacelike surfaces lying
 on rotational hypersurfaces in $\R^4_1$.
We shall consider the standard rotational hypersurface  with timelike axis and the rotational hypersurface with
spacelike axis.

Let $\{e_1, e_2, e_3, e_4\}$ be the standard orthonormal frame in
$\R^4_1$, i.e.  $e_1^2 =
e_2^2 = e_3^2 = 1, \, e_4^2 = -1$.
First we consider the standard rotational hypersurface with timelike axis.
Let $S^2(1)$ be a 2-dimensional sphere in the Euclidean space $\R^3 = \span
\{e_1, e_2, e_3\}$, centered at the origin $O$.
Let $f = f(u), \,\, g = g(u)$ be smooth functions, defined in an
interval $I \subset \R$, such that $\dot{f}^2(u) - \dot{g}^2(u) > 0, \,\, u \in I$,
where $\dot{f}(u)$ denotes the derivative $\ds{\frac{df(u)}{du}}$ and $\dot{g}(u)=\ds{\frac{dg(u)}{du}}$.
We assume that $f(u)>0, \,\, u \in I$. The standard rotational hypersurface $\mathcal{M}'$ in
$\R^4_1$, obtained by the rotation of the meridian curve $m: u
\rightarrow (f(u), g(u))$ about the $Oe_4$-axis,  is parameterized as follows:
$$\mathcal{M}': Z(u,w^1,w^2) = f(u)\,l(w^1,w^2) + g(u) \,e_4,$$
where $l(w^1,w^2)$ is the unit position vector of $S^2(1)$  in $\R^3$.
The hypersurface $\mathcal{M}'$ is a rotational hypersurface in $\R^4_1$ with timelike axis.

We consider a
smooth curve $c: l = l(v) = l(w^1(v),w^2(v)), \, v \in J, \,\, J \subset \R$  on
$S^2(1)$, parameterized by the arc-length, i.e. $\langle l'(v), l'(v) \rangle = 1$. Let
 $t(v) = l'(v)$  be the tangent vector field of $c$. Since $\langle t(v), t(v) \rangle = 1$,
 $\langle l(v), l(v) \rangle = 1$, and $\langle t(v), l(v) \rangle = 0$, there exists a unique (up to a sign)
 vector field $n(v)$ in $\R^3$, such that
$\{ l(v), t(v), n(v)\}$ is an orthonormal frame field. With respect to this
orthonormal frame field we have the following Frenet formulas of $c$ on $S^2(1)$:
\begin{equation} \label{E:Eq-5.1}
\begin{array}{l}
\vspace{2mm}
l' = t;\\
\vspace{2mm}
t' = \kappa \,n - l;\\
\vspace{2mm} n' = - \kappa \,t,
\end{array}
\end{equation}
 where $\kappa (v)= \langle t'(v), n(v) \rangle$ is the spherical curvature of $c$.

 Now we construct a surface $\mathcal{M}'_m$ in $\R^4_1$ in
the following way:
\begin{equation} \label{E:Eq-5.2}
\mathcal{M}'_m: z(u,v) = f(u) \, l(v) + g(u)\, e_4, \quad u \in I, \, v \in J.
\end{equation}
The surface $\mathcal{M}'_m$ lies on the rotational hypersurface $\mathcal{M}'$ in
$\R^4_1$. Since $\mathcal{M}'_m$ is a one-parameter system of meridians
of $\mathcal{M}'$, we call $\mathcal{M}'_m$ a
\textit{meridian surface on $\mathcal{M}'$}.

The tangent space of $\mathcal{M}'_m$ is spanned by the vector fields:
$$\begin{array}{l}
\vspace{2mm}
z_u = \dot{f} \,l + \dot{g}\,e_4;\\
\vspace{2mm} z_v = f\,t,
\end{array}$$
so, the coefficients of the first fundamental form of $\mathcal{M}'_m$
are $$E = \dot{f}^2(u) - \dot{g}^2(u) >0; \quad F = 0; \quad G = f^2(u) >0.$$
Hence, the first fundamental form is positive
definite, i.e. $\mathcal{M}'_m$ is a spacelike surface is $\R^4_1$.

Let us denote $x =\ds{\frac{z_u}{\sqrt{\dot{f}^2 - \dot{g}^2}}},\,\,y = \ds{\frac{z_v}{f}}$.
Then we have the following orthonormal tangent frame field of $\mathcal{M}'_m$:
$$x =\frac{\dot{f}(u)}{\sqrt{\dot{f}^2(u) - \dot{g}^2(u)}}\,\,l(v) + \frac{\dot{g}(u)}{\sqrt{\dot{f}^2(u) - \dot{g}^2(u)}}\, \, e_4; \qquad y =  t(v).$$
We consider the orthonormal normal frame field, defined by:
$$n_1 = n(v); \qquad n_2 = \frac{\dot{g}(u)}{\sqrt{\dot{f}^2(u) - \dot{g}^2(u)}}\,\,l(v) + \frac{\dot{f}(u)}{\sqrt{\dot{f}^2(u) - \dot{g}^2(u)}}\, \, e_4.$$
Thus we obtain a frame field $\{x,y, n_1,
n_2\}$ of $\mathcal{M}'_m$, such that $\langle n_1, n_1 \rangle =1$, $\langle n_2, n_2 \rangle =- 1$, $\langle n_1, n_2 \rangle =0$.

Taking into account
\eqref{E:Eq-5.1}, we calculate the second partial derivatives of $z(u,v)$:
$$\begin{array}{l}
\vspace{2mm}
z_{uu} = \ddot{f} \,l + \ddot{g}\,e_4;\\
\vspace{2mm}
z_{uv} = \dot{f}\,t;\\
\vspace{2mm} z_{vv} = f \kappa \,n - f\,l.
\end{array}$$
Then we get
$$\begin{array}{lll}
\vspace{2mm}
c_{11}^1 = \langle z_{uu}, n_1 \rangle = 0; & \qquad   c_{12}^1 = \langle z_{uv}, n_1 \rangle = 0; & \qquad  c_{22}^1 = \langle z_{vv}, n_1 \rangle = f \,\kappa;\\
\vspace{2mm}
c_{11}^2 = \langle z_{uu}, n_2 \rangle = \ds{\frac{\ddot{f} \dot{g} - \ddot{g} \dot{f} }{\sqrt{\dot{f}^2 - \dot{g}^2}}}; & \qquad  c_{12}^2 = \langle z_{uv}, n_2 \rangle = 0; & \qquad
c_{22}^2 = \langle z_{vv}, n_2 \rangle = \ds{- \frac{f \dot{g}}{\sqrt{\dot{f}^2 - \dot{g}^2}}}.
\end{array} $$
Hence,
\begin{equation}\label{E:Eq-5.3}
\begin{array}{l}
\vspace{2mm}
\sigma(x,x) =  \qquad \quad \; \ds{\frac{\dot{f} \ddot{g} - \dot{g}\ddot{f}}{(\dot{f}^2 - \dot{g}^2)^{\frac{3}{2}}}} \,\,n_2; \\
\vspace{2mm}
\sigma(x,y) = 0;  \\
\vspace{2mm}
 \sigma(y,y) = \ds{\frac{\kappa}{f}}\, \,n_1 + \ds{ \frac{\dot{g}}{f \sqrt{\dot{f}^2 - \dot{g}^2}}}\,\, n_2.
\end{array}
\end{equation}
Let us denote by $\kappa_m$ the curvature of the
meridian curve $m$, i.e. $\kappa_m (u)= \ds{\frac{\dot{f}(u) \ddot{g}(u) - \dot{g}(u) \ddot{f}(u)}{(\dot{f}^2(u) - \dot{g}^2(u))^{\frac{3}{2}}}}$.
Taking into account
\eqref{E:Eq-5.3}, we find the invariants $k$, $\varkappa$, and the Gauss curvature $K$ of $\mathcal{M}'_m$:
\begin{equation} \notag
k = - \frac{\kappa_m^2(u) \, \kappa^2(v)}{f^2(u)}; \qquad \varkappa = 0;
\qquad K = - \frac{\kappa_m (u)\, \dot{g}(u)}{f(u) \sqrt{\dot{f}^2(u) - \dot{g}^2(u)}}.
\end{equation}

The equality $\varkappa = 0$ implies that $\mathcal{M}'_m$ is a surface with
flat normal connection.

The mean curvature vector field $H$ is given by
\begin{equation}\label{E:Eq-5.5}
H = \frac{\kappa}{2f}\,\, n_1 + \frac{\kappa_m f \sqrt{\dot{f}^2 - \dot{g}^2} + \dot{g}}{2f \sqrt{\dot{f}^2 - \dot{g}^2}} \,\, n_2.
\end{equation}

Equality \eqref{E:Eq-5.5} implies that  $\mathcal{M}'_m$ is a minimal surface (the mean curvature vector field $H$ is zero) if and only if
$$\kappa = 0; \qquad \kappa_m f \sqrt{\dot{f}^2 - \dot{g}^2} + \dot{g} =0.$$

We shall exclude this case and consider the case when $\mathcal{M}'_m$ is marginally trapped, i.e. $H \neq 0$ and $\langle H,H \rangle =0$.
Using \eqref{E:Eq-5.5} we obtain the following result.

\begin{prop}\label{P:Meridian I}
The meridian surface  $\mathcal{M}'_m$, defined by \eqref{E:Eq-5.2},
 is marginally trapped if and only if
$$\kappa^2 = \ds{\frac{\left(\kappa_m f \sqrt{\dot{f}^2 - \dot{g}^2} + \dot{g}\right)^2}{\dot{f}^2 - \dot{g}^2}}, \qquad \kappa \neq 0.$$
\end{prop}

\vskip 5mm
In a similar way we shall consider meridian surfaces lying on the standard rotational hypersurface in $\R^4_1$
with spacelike axis.
Let $S^2_1(1)$ be the  timelike sphere in the Minkowski space $\R^3_1 = \span
\{e_2, e_3, e_4\}$, i.e. $S^2_1(1) = \{ V \in \R^3_1: \langle V, V \rangle = 1\}$.
$S^2_1(1)$ is a timelike surface in $\R^3_1$ known as the \emph{de Sitter space}.
Let $f = f(u), \,\, g = g(u)$ be smooth functions, defined in an
interval $I \subset \R$, such that $\dot{f}^2(u) + \dot{g}^2(u) >0$,  $f(u)>0, \,\, u \in I$.
We denote by $l(w^1,w^2)$  the unit position vector of $S^2_1(1)$  in $\R^3_1$ and
consider the rotational hypersurface $\mathcal{M}''$ in
$\R^4_1$, obtained by the rotation of the meridian curve $m: u
\rightarrow (f(u), g(u))$ about the $Oe_1$-axis. It is parameterized as follows:
$$\mathcal{M}'': Z(u,w^1,w^2) = f(u)\,l(w^1,w^2) + g(u) \,e_1.$$
The hypersurface $\mathcal{M}''$ is a rotational hypersurface in $\R^4_1$ with spacelike axis.

Now we consider a
smooth spacelike curve $c: l = l(v) = l(w^1(v),w^2(v)), \, v \in J, \,\, J \subset \R$  on
$S^2_1(1)$, parameterized by the arc-length, i.e. $\langle l'(v), l'(v) \rangle = 1$, and denote by
 $t(v) = l'(v)$ the tangent vector field of $c$. Since $\langle t(v), t(v) \rangle = 1$,
 $\langle l(v), l(v) \rangle = 1$, and $\langle t(v), l(v) \rangle = 0$, there exists a unique (up to a sign)
timelike vector field $n(v)$ in $\R^3_1$, such that $\langle n(v), n(v) \rangle = -1$,
$\langle n(v), t(v) \rangle = 0$, $\langle n(v), l(v) \rangle = 0$. In such a way we obtain
an orthonormal frame field
$\{ l(v), t(v), n(v)\}$ in $\R^3_1$. With respect to this
 frame field we have the following decompositions of the vector fields $l'(v)$, $t'(v)$, $n'(v)$:
\begin{equation} \notag
\begin{array}{l}
\vspace{2mm}
l' = t;\\
\vspace{2mm}
t' = - \kappa \,n - l;\\
\vspace{2mm} n' = - \kappa \,t,
\end{array}
\end{equation}
which can be considered as Frenet formulas of $c$ on $S^2_1(1)$.
The function $\kappa (v)= \langle t'(v), n(v) \rangle$ is the spherical curvature of $c$ on $S^2_1(1)$.

 We construct a surface $\mathcal{M}''_m$ in $\R^4_1$ in
the following way:
\begin{equation} \label{E:Eq-5.7}
\mathcal{M}''_m: z(u,v) = f(u) \, l(v) + g(u)\, e_1, \quad u \in I, \, v \in J.
\end{equation}
The surface $\mathcal{M}''_m$, defined by \eqref{E:Eq-5.7}, lies on the rotational hypersurface $\mathcal{M}''$ in
$\R^4_1$. Since $\mathcal{M}''_m$ is a one-parameter system of meridians
of $\mathcal{M}''$, we call $\mathcal{M}''_m$ a
\textit{meridian surface on $\mathcal{M}''$}.

The meridian surface $\mathcal{M}''_m$ is a spacelike surface is $\R^4_1$ with tangent vector fields
$$\begin{array}{l}
\vspace{2mm}
z_u = \dot{f} \,l + \dot{g}\,e_1;\\
\vspace{2mm} z_v = f\,t,
\end{array}$$
and  coefficients of the first fundamental form given by
$$E = \dot{f}^2(u) + \dot{g}^2(u) >0; \quad F = 0; \quad G = f^2(u) >0.$$

We consider the orthonormal tangent frame field
$x =\ds{\frac{z_u}{\sqrt{\dot{f}^2 + \dot{g}^2}}},\,\,y = \ds{\frac{z_v}{f}}$, i.e.
$$x =\frac{\dot{f}(u)}{\sqrt{\dot{f}^2(u) + \dot{g}^2(u)}}\,\,l(v) + \frac{\dot{g}(u)}{\sqrt{\dot{f}^2(u) + \dot{g}^2(u)}}\, \, e_1; \qquad y =  t(v),$$
and the orthonormal normal frame field, defined by:
$$n_1 = \frac{\dot{g}(u)}{\sqrt{\dot{f}^2(u) + \dot{g}^2(u)}}\,\,l(v) - \frac{\dot{f}(u)}{\sqrt{\dot{f}^2(u) + \dot{g}^2(u)}}\, \, e_1; \qquad  n_2 = n(v).$$
Thus we obtain a frame field $\{x,y, n_1, n_2\}$ of $\mathcal{M}''_m$,
such that $\langle n_1, n_1 \rangle =1$, $\langle n_2, n_2 \rangle =- 1$, $\langle n_1, n_2 \rangle =0$.

In the same way as in the previous example,  we obtain the
formulas
\begin{equation}\notag
\begin{array}{l}
\vspace{2mm}
\sigma(x,x) = \ds{\frac{\dot{g}\ddot{f} - \dot{f} \ddot{g}}{(\dot{f}^2 + \dot{g}^2)^{\frac{3}{2}}}} \,\,n_1; \\
\vspace{2mm}
\sigma(x,y) = 0;  \\
\vspace{2mm}
 \sigma(y,y) = \ds{- \frac{\dot{g}}{f \sqrt{\dot{f}^2 + \dot{g}^2}}}\,\, n_1 - \ds{\frac{\kappa}{f}}\, \,n_2.
\end{array}
\end{equation}
Now,  the curvature  $\kappa_m$ of the meridian curve $m$ is given
by $\kappa_m (u)= \ds{\frac{\dot{f}(u) \ddot{g}(u) - \dot{g}(u)
\ddot{f}(u)}{(\dot{f}^2(u) + \dot{g}^2(u))^{\frac{3}{2}}}}$. The
invariants $k$, $\varkappa$, and the Gauss curvature $K$ of
$\mathcal{M}''_m$ are expressed as follows:
\begin{equation} \notag
k = - \frac{\kappa_m^2(u) \, \kappa^2(v)}{f^2(u)}; \qquad
\varkappa = 0; \qquad K = \frac{\kappa_m (u)\, \dot{g}(u)}{f(u)
\sqrt{\dot{f}^2(u) + \dot{g}^2(u)}}.
\end{equation}

The equality $\varkappa = 0$ implies that $\mathcal{M}''_m$ is a surface with
flat normal connection.

The mean curvature vector field $H$ is given by
\begin{equation}\label{E:Eq-5.10}
H =  - \frac{\kappa_m f \sqrt{\dot{f}^2 + \dot{g}^2} + \dot{g}}{2f \sqrt{\dot{f}^2 + \dot{g}^2}} \,\, n_1 - \frac{\kappa}{2f}\,\, n_2.
\end{equation}

Obviously,   $\mathcal{M}''_m$ is a minimal surface (the mean curvature vector field $H$ is zero) if and only if
$$\kappa = 0; \qquad \kappa_m f \sqrt{\dot{f}^2 + \dot{g}^2} + \dot{g} =0.$$

We consider the case when $\mathcal{M}''_m$ is marginally trapped, i.e. $H \neq 0$ and $\langle H,H \rangle =0$.
Equality \eqref{E:Eq-5.10} implies the following result.

\begin{prop}\label{P:Meridian II}
The meridian surface  $\mathcal{M}''_m$, defined by \eqref{E:Eq-5.7},
 is marginally trapped if and only if
$$\kappa^2 = \ds{\frac{\left(\kappa_m f \sqrt{\dot{f}^2 + \dot{g}^2} + \dot{g}\right)^2}{\dot{f}^2 + \dot{g}^2}}, \qquad \kappa \neq 0.$$
\end{prop}

\vskip 5mm
At the end of this section we shall find all  marginally trapped meridian surfaces lying on the rotational hypersurfaces
$\mathcal{M}'$ or $\mathcal{M}''$.

Let $\mathcal{M}'_m$ be a marginally trapped meridian surface lying on the rotational hypersurface $\mathcal{M}'$.
According to Proposition \ref{P:Meridian I} the following equality holds:
\begin{equation}\notag
\kappa^2(v) = \ds{\frac{\left(\kappa_m(u) f(u) \sqrt{\dot{f}^2(u) - \dot{g}^2(u)} + \dot{g}(u)\right)^2}{\dot{f}^2(u) - \dot{g}^2(u)}},
\end{equation}
which imply
\begin{equation}  \label{E:Eq-5.11}
\begin{array}{l}
\vspace{2mm}
\kappa(v) = a = const, \qquad a \neq 0; \\
\vspace{2mm} \ds{\frac{\kappa_m(u) f(u) \sqrt{\dot{f}^2(u) - \dot{g}^2(u)} + \dot{g}(u)}{\sqrt{\dot{f}^2(u) - \dot{g}^2(u)}}} = \pm a.
\end{array}
\end{equation}

Without loss of generality we assume that the meridian curve $m$ is given by $$f = u; \qquad g = g(u).$$
Then $\dot{f}(u) = 1$; $\ddot{f}(u) = 0$; $\kappa_m(u) = \ds{\frac{\ddot{g}}{(1 - \dot{g}^2)^{\frac{3}{2}}}}$.
Hence, equation \eqref{E:Eq-5.11} takes the form:
\begin{equation}  \label{E:Eq-5.12}
u \ddot{g} + \dot{g} - (\dot{g})^3 = \pm a (1 - \dot{g}^2)^{\frac{3}{2}}.
\end{equation}

The meridian curves of all marginally trapped meridian surfaces lying on the rotational hypersurface $\mathcal{M}'$ are
determined by the solutions of differential equation \eqref{E:Eq-5.12}. After the change
$\ds{\frac{1}{1 - \dot{g}^2(u)} =  h^2(u) + 1}$ the above equation is transformed into
\begin{equation}  \label{E:Eq-5.13}
u  \dot{h}(u) + h(u)= \pm a \varepsilon, \qquad \varepsilon = \pm 1.
\end{equation}
The general solution of equation \eqref{E:Eq-5.13} is given by
\begin{equation} \notag
h(u)= \frac{\pm a \varepsilon u + c_1}{u}, \qquad c_1 = const.
\end{equation}
We set $c = \varepsilon c_1$ and get
\begin{equation}  \label{E:Eq-5.14}
\dot{g}(u) = \frac{\pm a  u + c}{\sqrt{(\pm a  u + c)^2 + u^2}}.
\end{equation}

Integrating \eqref{E:Eq-5.14} we obtain that all solutions of differential equation \eqref{E:Eq-5.12} are given by the formula
\begin{equation}  \label{E:Eq-5.15}
g(u) = \frac{\pm a}{a^2+1}\sqrt{(\pm a  u + c)^2 + u^2} +
\frac{c}{(a^2+1)^{\frac{3}{2}}} \ln \left(\sqrt{a^2+1}\, u \pm \frac{ac}{\sqrt{a^2+1}} + \sqrt{(\pm a  u + c)^2 + u^2} \right) + b,
\end{equation}
where $b$ and $c$ are constants, $c \neq 0$.

Thus we obtain the next result.

\begin{thm}\label{T:Meridian I}
The meridian surface  $\mathcal{M}'_m$
 is marginally trapped if and only if
 the curve $c$  on $S^2(1)$ has constant
spherical curvature $\kappa = a, \; a \neq 0$, and the
meridian $m$ is defined by
 \eqref{E:Eq-5.15}.
\end{thm}

Let $\mathcal{M}'_m$ be a marginally trapped meridian surface.
Then the meridian curve $m$ is determined by $f = u, \,\, g =g(u)$,
where the function $g(u)$ is defined by  \eqref{E:Eq-5.15}.
Hence, the curvature of $m$ is $\kappa_m = \ds{-\frac{c}{u^2}}$.
Formulas \eqref{E:Eq-5.3} take the form
\begin{equation} \label{E:Eq-5.16}
\begin{array}{l}
\vspace{2mm}
\sigma(x,x) =  \qquad \quad \; \ds{- \frac{c}{u^2}} \,\,n_2; \\
\vspace{2mm}
\sigma(x,y) = 0;  \\
\vspace{2mm}
 \sigma(y,y) = \ds{\frac{a}{u}}\, \,n_1 + \ds{ \frac{\pm a  u + c}{u^2}}\,\, n_2,
\end{array}
\end{equation}
which imply that the invariants $k$, $\varkappa$ and $K$ of $\mathcal{M}'_m$  are expressed as
\begin{equation}\notag
k = -\frac{a^2c^2}{u^6}; \qquad \varkappa = 0; \qquad K = \frac{c (\pm a  u + c)}{u^4}.
\end{equation}
The mean curvature vector field of $\mathcal{M}'_m$ is
$H = \ds{\frac{a}{2u}(n_1 \pm n_2)}$.

Note that $\{x,y,n_1,n_2\}$ is not the  geometric frame field of $\mathcal{M}'_m$ defined in Section \ref{S:Fundamental Thm}.
The geometric frame field  $\{X,Y,N_1,N_2\}$ of $\mathcal{M}'_m$ is determined by
\begin{equation}\notag
\begin{array}{ll}
\vspace{2mm}
X = \ds{\frac{1}{\sqrt{2}}(x+y)}; & \qquad  N_1 = \ds{\frac{a}{2u}(n_1 \pm n_2)}; \\
\vspace{2mm}
 Y = \ds{\frac{1}{\sqrt{2}}(x-y)}; & \qquad N_2 = \ds{\frac{u}{a}(- n_1 \pm n_2)}.
\end{array}
\end{equation}
Hence, using  \eqref{E:Eq-5.16} we obtain the formulas corresponding to the geometric frame field:
\begin{equation} \notag
\begin{array}{l}
\vspace{2mm}
\sigma(X,X) =  N_1; \\
\vspace{2mm}
\sigma(X,Y) = \ds{\frac{-au \mp c}{au}}\, \,N_1 \mp \ds{\frac{ac}{2 u^3}}\,\, N_2;  \\
\vspace{2mm}
 \sigma(Y,Y) = N_1.
\end{array}
\end{equation}
Thus we obtain the invariants $\nu$, $\lambda$, $\mu$ in the
Frenet-type derivative formulas of $\mathcal{M}'_m$:
\begin{equation}\notag
\nu = 0; \qquad \lambda = \ds{\frac{-au \mp c}{au}}; \qquad \mu =\ds{\mp \frac{ac}{2 u^3}}.
\end{equation}
The invariants $\gamma_1$, $\gamma_2$, $\beta_1$, $\beta_2$ are defined by
$$\gamma_1 = \langle\nabla'_XX, Y \rangle; \quad \gamma_2 = \langle\nabla'_YY, X \rangle; \quad
\beta_1 = - \langle\nabla'_XN_1, N_2 \rangle; \quad \beta_2 = -\langle\nabla'_YN_1, N_2 \rangle.$$
Calculating these scalar products we get
$$ \gamma_ 1 = \gamma_2 = \ds{-\frac{\sqrt{(\pm au + c)^2 + u^2}}{\sqrt{2}u^2}}; \qquad
\beta_1 = \beta_2 = \ds{\frac{\sqrt{(\pm au + c)^2 + u^2}}{\sqrt{2} u^4} \left(c(\pm au + c) - u^2\right)}.$$

Hence, the marginally trapped meridian surface $\mathcal{M}'_m$  has non-parallel mean curvature vector field,
since $\beta_1$ and $\beta_2$ are non-zero functions.

\vskip 4mm
In a similar way we  find all  marginally trapped meridian surfaces lying on the rotational hypersurface
$\mathcal{M}''$.

Let $\mathcal{M}''_m$ be a marginally trapped meridian surface lying on $\mathcal{M}''$.
Applying Proposition \ref{P:Meridian II} we obtain the following conditions
on the meridian curve $m$ and the curve $c$ lying on $S^2_1(1)$:
\begin{equation}  \notag
\begin{array}{l}
\vspace{2mm}
\kappa(v) = a = const, \qquad a \neq 0; \\
\vspace{2mm} \ds{\frac{\kappa_m(u) f(u) \sqrt{\dot{f}^2(u) + \dot{g}^2(u)} + \dot{g}(u)}{\sqrt{\dot{f}^2(u) + \dot{g}^2(u)}}} = \pm a.
\end{array}
\end{equation}

If the meridian curve $m$ is given by $f = u; \,\, g = g(u)$, then we get that
the function $g(u)$ is a solution of the following differential equation:
\begin{equation}\notag
u \ddot{g} + \dot{g} + (\dot{g})^3 = \pm a (1 + \dot{g}^2)^{\frac{3}{2}}.
\end{equation}
All solutions of this equation
are given by the formula
\begin{equation}  \label{E:Eq-5.17}
g(u) = \frac{\pm a}{1-a^2}\sqrt{u^2 - (\pm a  u + c)^2} +
\frac{c}{(1-a^2)^{\frac{3}{2}}} \ln \left(\sqrt{1-a^2}\, u \mp \frac{ac}{\sqrt{1-a^2}} + \sqrt{u^2 - (\pm a  u + c)^2} \right) + b,
\end{equation}
where $b$ and $c$ are constants, $c \neq 0$.

Thus we obtain the following result.

\begin{thm}\label{T:Meridian I}
The meridian surface  $\mathcal{M}''_m$
 is marginally trapped if and only if
 the curve $c$  on $S^2_1(1)$ has constant
spherical curvature $\kappa = a, \; a \neq 0$, and the
meridian $m$ is defined by
 \eqref{E:Eq-5.17}.
\end{thm}

The invariants of the marginally trapped meridian surface $\mathcal{M}''_m$ are expressed in a
similar way as the invariants of the marginally trapped meridian surface $\mathcal{M}'_m$.
The marginally trapped meridian surface $\mathcal{M}''_m$  has non-parallel mean curvature vector field,
since the invarianst $\beta_1$ and $\beta_2$ are non-zero.

\vskip 5mm
\textbf{Acknowledgements:} The second author is partially supported by "L.
Karavelov" Civil Engineering Higher School, Sofia, Bulgaria under
Contract No 10/2010.

\end{document}